\documentclass{ifacconf}

\usepackage{mypackages}

\newcommand{\beq}{\begin{equation}}
\newcommand{\eeq}{\end{equation}}
\newcommand{\beqa}{\begin{eqnarray}}
\newcommand{\eeqa}{\end{eqnarray}}
\newcommand{\beqan}{\begin{eqnarray*}}
\newcommand{\eeqan}{\end{eqnarray*}}

\newcommand{\rank}{\text{rank }}

\newcommand{\diag}{\mathop{\mathrm{diag}}}

\newcommand{\Acal}{{\cal A}}

\newcommand{\Ecal}{{\cal E}}
\newcommand{\Fcal}{{\cal F}}
\newcommand{\Gcal}{{\cal G}}

\newcommand{\Ncal}{{\cal N}}

\newcommand{\Vcal}{{\cal V}}

\newcommand{\Xcal}{{\cal X}}

\newcommand{\Zcal}{{\cal Z}}

\newcommand{\Gfrak}{\mathfrak{G}}

\newcommand{\Asf}{\sf{A}}

\newcommand{\Esf}{\sf{E}}

\newcommand{\bone}{\mathds{1}}

\renewcommand{\v}[1]{{\bm{#1}}}

\newcommand\T{{\mathpalette\raiseT\intercal}}
\newcommand\raiseT[2]{\raisebox{0.25ex}{$#1#2$}
}

\newcounter{l1}
\newcounter{l2}
\newcounter{l3}
\newcommand{\bdotlist}{\begin{list}{$\bullet$}{}}
\newcommand{\bboxlist}{\begin{list}{$\Box$}{}}
\newcommand{\bbboxlist}{\begin{list}{\raisebox{.005in}{{\tiny
$\blacksquare$ \ \ }}}{}}
\newcommand{\bdashlist}{\begin{list}{$-$}{} }
\newcommand{\blist}{\begin{list}{}{} }
\newcommand{\barablist}{\begin{list}{\arabic{l1}}{\usecounter{l1}}}
\newcommand{\balphlist}{\begin{list}{(\alph{l2})}{\usecounter{l2}}}
\newcommand{\bAlphlist}{\begin{list}{\Alph{l2}.}{\usecounter{l2}}}
\newcommand{\bdiamlist}{\begin{list}{$\diamond$}{}}
\newcommand{\bromalist}{\begin{list}{(\roman{l3})}{\usecounter{l3}}}

\newtheorem{theorem}{Theorem}
\newtheorem{lemma}{Lemma}
\newtheorem{definition}{Definition}

\newcommand{\proj}{\textrm{proj}}
\renewcommand{\hat}{\widehat}
\renewcommand{\bar}{\overline}
\renewcommand{\bone}{\mathds{1}}

\DeclareFontFamily{OT1}{pzc}{}
\DeclareFontShape{OT1}{pzc}{m}{it}{<-> s * [0.900] pzcmi7t}{}
\DeclareMathAlphabet{\mathpzc}{OT1}{pzc}{m}{it}

\begin{document}
\begin{frontmatter}
\title{On Privatizing Equilibrium Computation in Aggregate Games over Networks} 
\thanks[footnoteinfo]{This work is partially supported by a Siebel Energy Institute Grant and the Joan and Lalit Bahl Fellowship at Illinois.}
\author{Shripad Gade \qquad Anna Winnicki \qquad Subhonmesh Bose}
\address{Department of Electrical and Computer Engineering, University of Illinois at Urbana-Champaign, Urbana, IL 61801. \\ \{gade3,  annaw5, boses\}@illinois.edu}

\begin{abstract}               
We propose a distributed algorithm to compute an equilibrium in aggregate games where players communicate over a fixed undirected network. Our algorithm exploits correlated perturbation to obfuscate information shared over the network. We prove that our algorithm does not reveal private information of players to an honest-but-curious adversary who monitors several nodes in the network. In contrast with differential privacy based algorithms, our method does not sacrifice accuracy of equilibrium computation to provide privacy guarantees.
\end{abstract}

\end{frontmatter}

%!TEX root = root.tex

\section{Introduction}
\label{sec:intro}
Aggregate games are non-cooperative games in which a player's payoff or cost depends on her own actions and the sum-total of the actions taken by other players. In a Cournot oligopoly for example, firms compete to supply a product in a market with a price-responsive demand with a goal to maximize profit. A firm's profit depends on her production cost as well as the market price, where the latter only depends on the aggregate quantity of the product offered in the market by all firms. Aggregate games are widely studied in the literature, e.g., see \cite{novshek1985existence,jensen2010aggregative}. Multiple strategic interactions in practice admit an aggregate game model, e.g., Cournot competition models for wholesale electricity markets in \cite{willems2009cournot, cai2019role, cherukuri2019iterative}, supply function competition in general economies see \cite{jensen2010aggregative}, communication networks in \cite{teng2019application,koskie2005nash} and common agency games in \cite{martimort2011aggregate}. Aggregate games are often potential games and a pure-strategy Nash equilibrium can be guaranteed to exist. In this paper, we present an algorithm for networked players to compute such an equilibrium in a distributed fashion that maintains the privacy of players' cost structures.

Players in a networked game can only communicate with neighboring players in a communication graph. Distributed algorithms for computing Nash equilibrium in networked games have a rich literature, e.g., see \cite{koshal2016distributed,salehisadaghiani2018distributed,ye2017distributed,tatarenko2018accelerated,parise2015network}. The obvious difficulty in computing equilibrium strategy arises due to the inability of a player to observe the aggregate decision. Naturally distributed Nash computation proceeds via iterative estimation of the aggregate decision followed by local payoff maximization (or cost minimization) with a given aggregate estimate. \cite{koshal2016distributed,parise2015network} exploits consensus based averaging, \cite{koshal2016distributed,salehisadaghiani2018distributed} explore gossip based averaging, and \cite{tatarenko2018accelerated} employs gradient play along with acceleration for aggregate estimation over networks. 

%In this paper, we focus on distributed equilibrium computation in networked aggregate games. 
\vspace*{-0.05in}
\subsection{Our Contributions} 
\vspace{-0.05in}

%\vspace*{-0.1in}
%\begin{itemize}[leftmargin=*]
%\item 
Algorithms for equilibrium computation were not designed with privacy in mind. We show in Section \ref{Sec:NoPrivacy}, that an honest-but-curious adversary can compromise a few nodes in the network and observe the sequence of estimates to infer other players' payoff or cost structures for the algorithm in \cite{koshal2016distributed}. In other words, information that allows distributed equilibrium computation can leak players' sensitive private information to adversaries. 
%Our analysis extends to show privacy breach in other distributed Nash equilibrium computation algorithms. 
%While we focus on the algorithm in \cite{koshal2016distributed}, our analysis can be extended to prove that algorithms in \cite{salehisadaghiani2018distributed,ye2017distributed,tatarenko2018accelerated,parise2015network} are not private either.

Distributed equilibrium computation algorithms require aggregate estimates to update their own actions. Our proposed algorithm obfuscates local aggregate estimates before sharing them with neighbors. The obfuscation step involves players adding correlated perturbations to each outgoing aggregate estimate. The perturbations are designed such that they add to zero for each player. The received perturbed aggregate estimates are averaged by each player and used for updating strategy using local projected gradient descent. 

%\vspace{0.05in}
%\item 
Our main result (Theorem \ref{Th:Main}) reveals that obfuscation via correlated perturbations prevents an adversary from accurately learning cost structures provided the network satisfies appropriate connectivity conditions. Players converge to exact Nash equilibrium asymptotically. In other words, we simultaneously achieve both privacy and accuracy in distributed Nash computation in aggregate games. This is in sharp contrast to differentially private algorithms where trade-offs between accuracy and privacy guarantee are fundamental, e.g., see \cite{7431982}. 

%\vspace{0.05in}
%\item  
Simulations in Section~\ref{Sec:Numerics} validate our results and corroborate our intuition that obfuscation slows down but does not impede the convergence of the algorithm. 
%\end{itemize}

%!TEX root = root.tex

\section{Equilibrium  Computation in Aggregate Games and the Lack of Privacy}
\label{sec:Model}
We begin by introducing a networked aggregate game. We then present an adversary model and show that prior distributed equilibrium computation algorithms leak private information of players. This exposition motivates the development of privacy-preserving algorithms for equilibrium computation in the next section.

\subsection{The Networked Aggregate Game Model} \label{Sec:Problem}
Consider a game with $N$ players that can communicate over a fixed undirected network with reliable lossless links. Model this communication network by graph  $\Gfrak(\Vcal, \Ecal)$, where each node in $\Vcal := \{1, \ldots, N\}$ denotes a player. Two players $i$ and $j$ can communicate with each other if and only if they share an edge in $\Ecal$, denoted as $(i,j)\in\Ecal$. Call $\Ncal_i$ the set of neighbors of node $i$ and $i\in\Ncal_i$ by definition.

Player $i$ can take actions in a convex compact set $\Xcal_i \subseteq \Rset^d$, where $\Rset$ denotes the set of real numbers. Define $\bar{\Xcal}$ as the Minkowski (set) sum of $\Xcal_i$'s and
$$\bar{x} \coloneqq \sum_{j=1}^N x_j$$
as the aggregate action of all players. For convenience, define $\bar{x}_{-i} \coloneqq \sum_{j \neq i}x_i$. We assume that $\cap_{i=1}^N \Xcal_i$ is non-empty. 
For an action profile $(x_1, \ldots, x_N)$, player $i$ incurs a cost that takes the form $f_i(x_i,\bar{x}) \coloneqq f_i \left(x_i, x_i + \bar{x}_{-i} \right)$. This defines an  aggregate game in that the actions of other players affect player $i$ only through the sum of actions of all players, $\bar{x}$. 

Each player $i\in\Vcal$ thus seeks to solve
\begin{alignat}{2}
\begin{aligned}
& \text{minimize}  && f_i(x_i,x_i+\bar{x}_{-i}), \\
& \text{subject to}  && x_i \in \Xcal_i.
\end{aligned}
\label{Eq:opt.i}
\end{alignat}

For each $i \in \Vcal$, assume that $f_i(x_i,y)$ is continuously differentiable in $(x_i,y)$ over a domain that contains $\Xcal_i \times \bar{\Xcal}$. Furthermore, for each $i\in\Vcal$, let $x_i \mapsto f_i(x_i, \bar{x})$ be convex over $\Xcal_i$ and the gradient $\nabla_{x_i} f_i$ be uniformly $\bar{L}$-Lipschitz, i.e.,  $\exists \ \bar{L} > 0$ such that,
\begin{align}
    \|\nabla_{x_i} f_i(x_i, u) - \nabla_{x_i} f_i(x_i, u')\| \leq \bar{L} \|u - u'\|, 
\end{align}
for all $u, u'$ in $\bar{\Xcal}$, $x_i$ in $\Xcal_i$.
Throughout, $\| \cdot \|$ stands  for  the $\ell_2$-norm of its argument. 
Define $\Xcal :=\times_{i=1}^N{\Xcal_i}$ and the gradient map 
\begin{align}
\phi(x) := \begin{pmatrix}
\nabla_{x_1} f_1(x_1,\bar{x}) \\
%\nabla_{x_2} f_2(x_2,\bar{x}) \\
\vdots \\
\nabla_{x_N} f_N(x_N,\bar{x}) \label{Eq:DefinePhi}
\end{pmatrix} 
\end{align}
for $x := (x_1^\T, x_2^\T, \ldots,x_N^\T)^\T \in \Xcal$.
Assume throughout that $\phi$ is strictly monotone over $\Xcal$, i.e., 
\begin{align}
\left[\phi(x) -\phi(x')\right]^\T(x-x') > 0, \label{Eq:StrictMonotoneGradMap}
\end{align}
for all $x,x' \in \Xcal$ and $x\neq x'$.
Denote this game in the sequel by ${\sf game}( \Gfrak, \{f_i, \Xcal_i\}_{i \in \Vcal} )$.

To provide a concrete example, consider the well-studied Nash-Cournot game (see \cite{fudenberg1991game}) among $N$ suppliers competing to offer into a market for a single commodity where the price $p$ varies with demand $D$ as $p(D) := a - b D$. Supplier $i$ offers to produce $x_i$ amount of goods within its production capability modeled as $\Xcal_i \subseteq \Rset_+$. Here $\Rset_+$ denotes the set of nonnegative real numbers. To produce $x_i$, supplier $i$ incurs a cost of $c_i(x_i)$, where $c_i$ is increasing, convex and differentiable. Each supplier seeks to maximize her profit, or equivalently, minimize her loss. The loss of supplier $i$ is 
$$ f_i(x_i,\bar{x})= c_i(x_i) - x_i p(\bar{x}) = c_i(x_i) - x_i(a - b\bar{x}).$$

\subsection{Equilibrium Definition and Existence}

An action profile $(x_1^*, \ldots 
x_N^*)$ defines a Nash equilibrium of ${\sf game}( \Gfrak, \{f_i, \Xcal_i\}_{i \in \Vcal} )$ in pure strategies, if
$$ 
f_i\left(x_i^*, x_i^* + \bar{x}_{-i}^* \right) \leq 
f_i\left(x_i, x_i + \bar{x}_{-i}^* \right),
$$
for all $x_i \in \Xcal_i$ and $i \in \Vcal$.

The networked aggregate game, as described above, always admits a unique pure strategy Nash equilibrium. See Theorem~2.2.3 in \cite{facchinei2007finite} for details. Given that an equilibrium always exists, prior literature has studied distributed algorithms for players to compute such an equilibrium. 
%in \cite{salehisadaghiani2018distributed,ye2017distributed,tatarenko2018accelerated,parise2015network}

\subsection{Prior Algorithms for Distributed Nash Computation } 
\label{Sec:PrivacyBreachAnalysis}

We now describe the distributed algorithm in \cite{koshal2016distributed} for equilibrium computation of ${\sf game}( \Gfrak, \{f_i, \Xcal_i\}_{i \in \Vcal} )$. In Section \ref{Sec:NoPrivacy}, we demonstrate that adversarial players can infer private information about cost structures $f_i$'s from observing a subset of the variables during equilibrium computation using that algorithm. While we only study the algorithm in \cite{koshal2016distributed}, our analysis can be extended to those presented in  \cite{salehisadaghiani2018distributed,ye2017distributed,tatarenko2018accelerated,parise2015network}.

Recall that players in ${\sf game}( \Gfrak, \{f_i, \Xcal_i\}_{i \in \Vcal} )$ do not have access to the aggregate decision. To allow equilibrium computation, let players at iteration $k$ maintain estimates of the aggregate decision $\bar{x}$ as $v^k_1, \ldots, v^k_N$, initialized as, $v^0_i = x^0_i$ for each player $i$. At discrete time steps $k \geq 0$, each player transmits her own estimate of the aggregate decision to its neighbors and updates her own action as,
\begin{subequations}
\begin{align}
\hat{v}^k_i &= \sum_{j=1}^N W_{ij} v^k_j,  \label{Eq:DAlgo1a}\\
x^{k+1}_i &= \proj_{\Xcal_i} \left[x^k_i - \alpha^k \nabla_{x_i} f_i(x^k_i, N \hat{v}^k_i) \right], \label{Eq:DAlgo1b}\\
v^{k+1}_i &= \hat{v}^k_i + x^{k+1}_i - x^k_i. \label{Eq:DAlgo1c}
\end{align}
\label{Eq:DAlgo1}%
\end{subequations}%%
Here, $\proj_{\Xcal_i}$ stands for projection on $\Xcal_i$, and $\alpha^k$ is a common learning rate of all players.

The algorithm has three steps. First, player $i$ computes a weighted average of the estimates of the aggregate received from its neighbors in \eqref{Eq:DAlgo1a}, where $W$ is a symmetric doubly-stochastic weighting matrix. The sparsity pattern of the matrix follows that of graph $\Gfrak$, i.e., 
$$ W_{ij} \neq 0 \iff (i, j) \in \Ecal.$$ 
Second, player $i$ performs a projected gradient update in \eqref{Eq:DAlgo1b} utilizing the weighted average of local aggregate decision $\hat{v}^k_i$ in lieu of the true aggregate decision $\bar{x}$. Finally, she updates her own estimate of aggregate average in \eqref{Eq:DAlgo1c} based on her  local decision $x_i^k$ and its update $x_i^{k+1}$.

\subsection{Adversary Model and Privacy Definition}  \label{Sec:AdvModelandPrivacyDef}
Consider an adversary ${\Asf}$ that compromises the players in $\Acal \subseteq \Vcal$. ${\Asf}$ is equipped with unbounded storage and computational capabilities, and has access to all information stored, processed locally and communicated to any compromised players at all times. We define adversary model using the information available to ${\Asf}$.
\begin{enumerate}[label=(\Alph*),leftmargin=*]
    \item  For a compromised node $i \in\Acal$, ${\Asf}$ knows all local information $f_i$, $x_i^k$, $v_i^k$, $\hat{v}_i^k$ and information received from neighbors of $i$ i.e., $v_j^k$ for $j \in \Ncal_i$ at each $k \geq 0$. 
    \item ${\Asf}$ knows the algorithm for equilibrium computation and its parameters $\{\alpha^k\}$ and $W$. 
    \item ${\Asf}$ observes aggregate decision $\bar{x}^k$ at each $k$.
\end{enumerate}

%while our algorithm guarantees privacy  against a strong adversary. 
% We purposefully bolster adversarial strength and perform privacy analysis in the worst case scenario.

%We use the notation $\kappa(\Gfrak)$ for vertex connectivity of graph $\Gfrak$. 

What does ${\Asf}$ seek to infer? The dependency of a player's cost on her own actions encodes private information. In the Cournot competition example, this dependency is precisely supplier $i$'s production cost -- information that is business sensitive. ${\Asf}$ seeks to exploit information sequence observed from compromised players to infer private information of other players. Intuitively, privacy implies inability of ${\Asf}$ to infer private cost functions. 
%and show that Algorithm~\ref{Eq:DAlgo1} is not private, implying ${\Asf}$ can learn the production costs of all players in a Cournot competition (Section \ref{Sec:NoPrivacy}). 

\begin{figure}[h]
\centering
\includegraphics[width=0.4\textwidth]{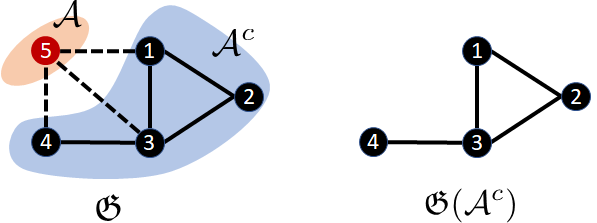}
\caption{Illustration of $\Gfrak$ and $\Gfrak(\Acal^c)$. Here, $\Acal = \{ 5 \}$ and $\Acal^c = \{ 1,2,3,4 \}$. }
\label{fig:graphs}
\end{figure}

Denote the set of non-adversarial nodes by $\Acal^c \coloneqq \Vcal \setminus \Acal$. Call $\Gfrak (\Acal^c)$ the restriction of $\Gfrak$ to $\Acal^c$ obtained by deleting the adversarial nodes. See Figure \ref{fig:graphs} for an illustration. For this example, $\Asf$ monitors all variables and parameters pertaining to player 5, but seeks to infer the functions $f_1, \ldots, f_4$. 

Let $\Pi$ denote the set of all permutations over all non-adversarial nodes in $\Acal^c$.
%, that maps each non-adversarial node to any other non-adversarial nodes, i.e. $\pi(i) = j$, if $i, j \notin \Acal$ ($i$ could be equal to $j$). We have $(N-\tau)!$ such possible permutations provided that $\tau$ out of $N$ nodes are corrupted by ${\Asf}$. 
%We use these permutations to construct set $\Fcal$, defined as the set of all network games obtained by permuting the cost function and decision set assignment among non-adversarial nodes. Mathematically,
Define the collection of games
\begin{align*} %\label{Eq:AggInvariantProblem}
	\Fcal \coloneqq \Big\{ {\sf game}( \Gfrak, \{f_{\pi(i)}, \Xcal_{\pi(i)}\}_{i \in \Vcal} ) \; \Big| \;\pi \in \Pi \Big\}.
\end{align*}
Thus, $\Fcal$ comprises the games where the cost functions and strategy sets of non-adversarial players are permuted. All games in $\Fcal$ have the same aggregate strategy $\bar{x}^*$ at Nash equilibrium. Next, we utilize $\Fcal$ to define privacy.

\vspace{0.1in}
\begin{definition}[Privacy] \label{Def:Privacy}
Consider a distributed algorithm to compute the Nash equilibrium of ${\sf game}( \Gfrak, \{f_i, \Xcal_i\}_{i \in \Vcal} )$. If execution observed by adversary ${\Asf}$ is consistent with all games in $\Fcal$, then the algorithm is private.
\end{definition}

We define privacy as the inability of ${\Asf}$ to distinguish between games in $\Fcal$.
%, i.e., if ${\Asf}$ potentially observing exact same execution for any game in $\Fcal$. %As a result we call $\Fcal$ as the set of aggregate invariant games.
Even if ${\Asf}$ knew all possible costs exactly--which is a tall order--our privacy definition implies that ${\Asf}$ cannot associate such costs to specific players.

\subsection{Privacy Breach in Algorithm \eqref{Eq:DAlgo1}} \label{Sec:NoPrivacy}

Consider a Cournot competition among 5 players connected according to $\Gfrak$ in Figure \ref{fig:graphs}, where ${\Asf}$ has compromised player $5$. Assume that the equilibrium of the game lies in the interior of each player's strategy set. Recall that $\Asf$ stores observed information at each $k$ and processes it to infer private cost information $c_i(x_i)$. We argue how $\Asf$ can compute cost functions $c_1(\cdot), \ldots, c_4(\cdot)$ up to a constant.

We first show privacy breach for player 4. ${\Asf}$ observes $\{v^k_1, v^k_3, v^k_4,v^k_5\}$ at each $k \geq 0$. $\Asf$ uses $v^k_3$, $v^k_4$, $v^k_5$ and $W$ to compute $\hat{v}^k_4$ using \eqref{Eq:DAlgo1a}. Moreover, ${\Asf}$ uses \eqref{Eq:DAlgo1c} to compute, 
$$x^{k+1}_4 - x^k_4 = v^{k+1}_4 - \hat{v}^k_4.$$

For large enough $k$, the step-size $\alpha^k$ is small enough to ensure,
   {\small
    $$\proj_{\Xcal_i} \left[x^k_i - \alpha^k \nabla_{x_i} f_i(x^k_i, N \hat{v}^k_i) \right] \approx x^k_i - \alpha^k \nabla_{x_i} f_i(x^k_i, N \hat{v}^k_i).$$}%
At such large $k$, ${\Asf}$ uses 
%$x^{k+1}_4 - x^k_4 = \alpha^k\nabla_{x_4} f_4(x^k_4,N\hat{v}^k_4)$, 
\eqref{Eq:DAlgo1b} along with $(x^{k+1}_4 - x^k_4)$ and $\alpha^k$ to calculate $\nabla_{x_4} f_4(x^k_4,N\hat{v}^k_4)$. 

${\Asf}$ uses information about strucutre of loss function i.e. $f_4(x_4,\bar{x}) = c_4(x_4) - x_4(a-b\bar{x})$, along with $\nabla_{x_4} f_4(x^k_4,N\hat{v}^k_4)$, $\hat{v}^k_4$, $\bar{x}^k$ and game parameters $a,b$ to learn $c'_4(x^k_4)$. Several observations of $(x^k_4,c'_4(x^k_4))$ allows ${\Asf}$ to learn the private cost $c_4$ upto a constant.

We showed that privacy breach for player 4, the same analysis can be used for players 1, 2 and 3 with an additional step. ${\Asf}$ observes $\bar{x}^k$, which tracks $\frac{1}{N}\sum_i v^k_i$ (Lemma 2 in \cite{koshal2016distributed}). ${\Asf}$ computes $$v^k_2 = N\bar{x}^k - (v^k_1+v^k_3+v^k_4+v^k_5).$$Since $\{v^k_2\}$ is available for each $k\geq0$, ${\Asf}$ uses same process as above to show privacy breach for players 1, 2 and 3.

For algorithm \eqref{Eq:DAlgo1}, $\Asf$ uncovers all private cost functions $c_i(\cdot)$ for an example aggregate game. Next, we design an algorithm that protects privacy of players' private information in the sense of Definition \ref{Def:Privacy} against  ${\Asf}$.

%!TEX root = root.tex

\section{Our Algorithm and Its Properties} 
\label{Sec:algProp}
We  propose and analyze Algorithm \ref{Algo:PrivDNEComp} that computes Nash equilibrium of ${\sf game}( \Gfrak, \{f_i, \Xcal_i\}_{i \in \Vcal} )$ in a distributed fashion. The main result (Theorem \ref{Th:Main}) shows that the algorithm asymptotically converges to the equilibrium. Attempts by ${\Asf}$ to recover each player's cost structure, however, remain unsuccessful.

The key idea behind our design is the injection of correlated noise perturbations in the exchange of local estimates of the aggregate decision. Different neighbors of player $i$ receive different estimates of the aggregate decision. The perturbations added by any player $i$ add to zero. While $\Asf$ may still infer the true aggregate decision, the protocol does not allow him to correctly infer the players' iterates or the gradients of their costs with respect to their own actions. Our assumption on network connectivity requires $\Gcal(\Acal^c)$ be connected and not be bipartite. Under these conditions $\Asf$ cannot monitor all outgoing communication channels from any player. We further show that one can design noises in a way that $\Asf$'s observations are consistent with all games in $\Fcal$, making it impossible for him to uncover cost for any specific player.

Throughout, assume that $W$ is a doubly stochastic that follows the sparsity pattern of $\Gfrak$. Further, assume that all non-diagonal, non-zero entries of $W$ are identically $\delta < \frac{1}{N-1}$.

At each time $k$, player $i$ generates correlated random numbers $\{r^k_{ij}\}$ satisfying $r^k_{ii}=0$ and $\sum_{j \in \mathcal{N}_i} r^k_{ij} = 0$. 
Player $i$ then adds $\alpha^k  r^k_{ij}$ to $v^k_{i}$ to generate $v^k_{ij}$, the estimate sent by player $i$ to player $j$, according to \eqref{Eq:PrivateAlgo1}. Let $\v{r}$ denote the collection of $r$'s for all players across time. Call $\v{r}$ the obfuscation sequence.

Each node $i$ computes weighted average of received aggregate estimates $v^k_{ji}$ to construct its own estimate aggregate decision $N\hat{v}^k_i$, following \eqref{Eq:PrivateAlgo2}. Players perform projected gradient descent using local decision estimate $x^k_i$, gradient of cost function $\nabla_{x_i} f_i(x^k_i,N\hat{v}^k_i)$, and non-summable, square-summable step size $\alpha^k$ (see \eqref{Eq:Stepsize}) to arrive at an improved local decision estimate $x^{k+1}_i$ using \eqref{Eq:PrivateAlgo3}. Players then update their local aggregate estimate using the change in local decision estimate $x^{k+1}_i - x^k_i$ per \eqref{Eq:PrivateAlgo4}.
\begin{algorithm}[!t]
\caption{Private Distributed Nash Computation}\label{Algo:PrivDNEComp}
\begin{algorithmic}[1]
\Statex \hspace*{-0.25in} \textbf{Input: }Player $i$ knows $f_i(x_i,\bar{x})$, $\mathcal{X}_i$, and $\delta$. Consider a non-increasing non-negative sequence $\v{\alpha}$ that satisfies 
\hspace*{-0.25in}\begin{align}
\hspace*{-0.25in}
%\alpha^k \geq 0, \ \alpha^k \geq \alpha^{k+1}, 
\sum_{k=1}^\infty \alpha^k = \infty \text{ and } \sum_{k=1}^\infty [\alpha^k]^2 < \infty. \label{Eq:Stepsize}
\end{align}
%\Statex \hspace*{-0.2in} \textbf{Result: }Iterates $\{x^k_i\}$ converge to NE denoted by $\{x^*_i\}$
%\Statex
\Statex \hspace*{-0.2in}\textbf{Initialize:} For $i\in\Vcal$, $v^0_i = x^0_i = \mathpzc{x} \in\cap_i \Xcal_i$.
\Statex
\Statex \hspace*{-0.25in} For $k\geq 0$, players $i\in\Vcal$ execute in parallel:
\Statex
\State Construct $|\mathcal{N}_i|$ random numbers $\{r^k_{ij}\}$, satisfying 
\begin{align}
r_{ii}^k = 0 \text{ and } \sum_{j \in \mathcal{N}_i} r^k_{ij} = 0.\label{Eq:PrivateAlgo0}   
%  \; \text{ and } \; |r^k_{ij} |\leq \Delta. 
\end{align}
\State Send obfuscated aggregate estimates $v^k_{ij}$ to $j \in \Ncal_i$, where
\begin{align}
    v^k_{ij} = v^k_i + \alpha^k r^k_{ij}.  \label{Eq:PrivateAlgo1}
\end{align}
\State Compute weighted average of received estimates $v^k_{ji}$ as
\begin{align}
\hat{v}^k_i = \sum_{j=1}^N W_{ij} v^k_{ji}. \label{Eq:PrivateAlgo2}
\end{align}
\State Perform a projected gradient descent step as
\begin{align}
x^{k+1}_i = \proj_{\mathcal{X}_i} [x^k_i - \alpha^k \nabla_{x_i} f_i(x^k_i, N \hat{v}^k_i)]. \label{Eq:PrivateAlgo3}
\end{align}
\State Update local aggregate estimate as
\begin{align}
v^{k+1}_i = \hat{v}^k_i + x^{k+1}_i - x^k_i. \label{Eq:PrivateAlgo4}
\end{align}
\end{algorithmic}
\end{algorithm}
The properties of our algorithm are summarized in the next result. The proof is included in Section \ref{Sec:proofs}. 

\vspace{0.075in}
\begin{theorem}
Consider a networked aggregate game defined as ${\sf game} ( \Gfrak, \{f_i, \Xcal_i\}_{i \in \Vcal})$. If $\Gfrak(\Acal^c)$ is connected and not bipartite, then Algorithm \ref{Algo:PrivDNEComp} is private. Moreover, if the obfuscation sequence is bounded, then Algorithm \ref{Algo:PrivDNEComp} asymptotically converges to a Nash equilibrium of the game.
\label{Th:Main}
\end{theorem}
\vspace{0.05in}

The convergence properties largely mimic that of distributed descent algorithms for  equilibrium computation. The locally balanced and bounded nature of the designed noise together with decaying step-sizes ultimately drown the effect of the noise. Computing balanced yet bounded perturbations can be achieved using secure multiparty computation protocols described in \cite{gade2016private,gade2018acc,abbe2012privacy}.
Our assumption on $\Gfrak(\Acal^c)$
is such that given two games $F, \tilde{F}$ from $\Fcal$ and an obfuscation sequence $\v{r}$, we are able to design a different obfuscation sequence $\v{\tilde{r}}$, such that the execution of $F$ perturbed with $\v{r}$ generates identical observables as $\tilde{F}$ perturbed with $\v{\tilde{r}}$. The connectivity among non-adversarial players in $\Acal^c$ is key to the success of our algorithm design. Convergence speed depends on the size of the perturbations. We investigate this link experimentally in Section~\ref{Sec:Numerics}, but leave analytical characterization of this relationship for future work. In what follows, we compare our algorithm and its properties to other protocols for privacy preservation.

%following adversarial model,

%\noindent \textbf{Privacy Results: }We motivated private games in Sec.~\ref{Sec:Problem}, where, we are interested in distributed computation of Nash equilibrium, although we aim to hide the cost function $f_i(x_i,\bar{x})$ from adversarial nodes.    
%\begin{theorem}[Privacy] \label{Th:PrivacyResult}
%Any set of $\tau$ adversarial nodes cannot accurately infer $c_i(x_i)$ if and only if the vertex connectivity of graph is greater than $\tau + 1$. 
%\end{theorem}
%The privacy result (Theorem~\ref{Th:PrivacyResult}) is a consequence of the fact that each node shares a perturbed aggregate estimate. And the perturbed aggregates appear to have been generated by a node with perturbed cost function. Adversary can infer this perturbed cost function but not the original cost function. The graph condition (vertex connectivity $\geq \tau+1$) is both necessary and sufficient. Intuitively, deletion of $\tau$ compromised nodes should not partition the non-adversarial nodes into disconnected components, as that may lead to adversary learning about the effective cost function of the components.  

%In what follows, we compare our result to existing privacy preserving computation in optimization and game theoretic settings.

\subsubsection{Comparison with Differentially Private Algorithms:}Differentially private algorithms for computing Nash equilibrium of potential games have been studied in  \cite{dong2015differential,cummings2015privacy}. The algorithm in \cite{dong2015differential} executes a differentially private distributed mirror-descent algorithm to optimize the potential function. 
% followed by using mirror-descent algorithm for computing an optimum. 
Experiments reveal that a trade-off arises between accuracy and privacy parameters, i.e., the more privacy one seeks, the less accurate the final output of the algorithm becomes. Such a tradeoff is a hallmark of differentially private algorithms, e.g., see \cite{7431982}.  Our algorithm on the other hand does not suffer from that limitation. Notice that our definition of privacy is binary in nature. That is, an algorithm for equilibrium computation can either be private or non-private. We aim to explore properties of our algorithmic architecture with notions of privacy that allow for a degree of privacy and compare them with differentially private algorithms.

%\begin{figure}[t]
%    \centering
%    \includegraphics[width=0.5\textwidth]{Fig2.PNG}
%    \caption{Graph with $\tau=1$ adversary in red and $\Hcal$ edges in Green. (Left) $\kappa(\Gcal)=3$. (Right) $\kappa(\Gcal)=2$.}
%    \label{Fig:Proof1}
%\end{figure} 

\subsubsection{Comparison to Cryptographic Methods:} Authors in \cite{lu2015game} use secure multiparty computation to compute Nash equilibrium. Such an approach guarantees privacy in an information theoretic sense. This protocol provides privacy guarantees along with accuracy, similar to our algorithmic framework. However, cryptographic protocols are typically computationally expensive for large problems (see Section V in \cite{zhang2019admm}), and are often difficult to implement in distributed settings. 

% Alternatively, we can use cryptographic protocols to privately build an aggregate estimate followed by local projected gradient descent step to get information theoretic privacy along with accuracy. However, as discussed above this approach incurs large computation and communication cost. 

\subsubsection{Comparison to Private Distributed Optimization:} Our earlier work in \cite{gade2018acc} has motivated the design of Algorithm \ref{Algo:PrivDNEComp}. While our prior work seeks privacy-preserving distributed protocols to cooperatively solve optimization problems, the current paper focuses on non-cooperative games. Protocols in \cite{gade2018acc} advocate use of perturbations that cancel over the network. Such a design is not appropriate for networked games for two reasons. First, players must agree on noise design, a premise that requires cooperation. Second, perturbing local functions $f_i$'s, even if the changes cancel in aggregate, can alter the equilibrium of the game.
%There is a crucial difference between the protocol in \cite{gade2018acc} and Algorithm \ref{Algo:PrivDNEComp}. The former advocates the use of perturbations  that zero out over the network as opposed to requiring them to be locally balanced.

\subsubsection{Privacy in Client-Server architecture:} This work considers players communicating over a peer-to-peer network. However, engineered distributed systems often have a client-server architecture. 
%Independent System Operators (ISO's) act as a server or a central coordinating entity and the generators act as clients. Economic dispatch over power grid is an aggregate game, where, generators are players and decide the quantity of power to generate. 
Presence of a central server entity allows for easy aggregate computation. However, privacy is sacrificed if the parameter server is adversarial. We have investigated privacy preservation for distribution optimization in this architecture in \cite{gade2018cdc}, where, we use multiple central servers instead of one, a subset of which can be adversarial. We believe our algorithm design and analysis in \cite{gade2018cdc} can be extended to deal with private equilibrium computation for aggregate games in client-server framework.
%Multiple parameter server architecture from  can be leveraged to privately aggregate local decisions at the server and improve privacy while ensuring fast NE computation (using access to aggregate decision).

%\subsubsection{Privacy and Convergence Rate Trade-off:}\hspace*{-0.1in}Privacy requires that noise bound $\Delta$ be large enough, to ensure that, any game in $\Fcal$ may lead to the same observations by ${\Asf}$. If $\Delta$ is smaller than required bound then some games in $\Fcal$ may lead to different observations by ${\Asf}$, and allow ${\Asf}$ to identify some games in $\Fcal$ as possible and some games as impossible. This reduces the size of games that lead to same observed execution by adversary from $\Fcal$ to a subset of $\Fcal$, and weakens privacy. Hence, there is a privacy and convergence rate trade-off, although after $\Delta$ is large enough, the privacy does not increase, albeit the convergence rate suffers. 

%!TEX root = root.tex

\section{A Numerical Experiment}
\label{Sec:Numerics}

\begin{figure}[!t]
    \centering
    \includegraphics[width=0.40\textwidth]{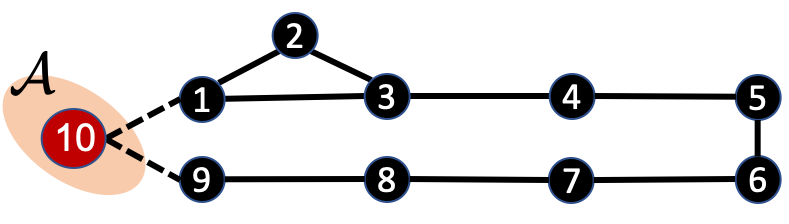}
    \caption{Communication network for Cournot network example on $N=10$ players.}
    \label{fig:graph_example}
    \centering
    \includegraphics[width=0.5\textwidth]{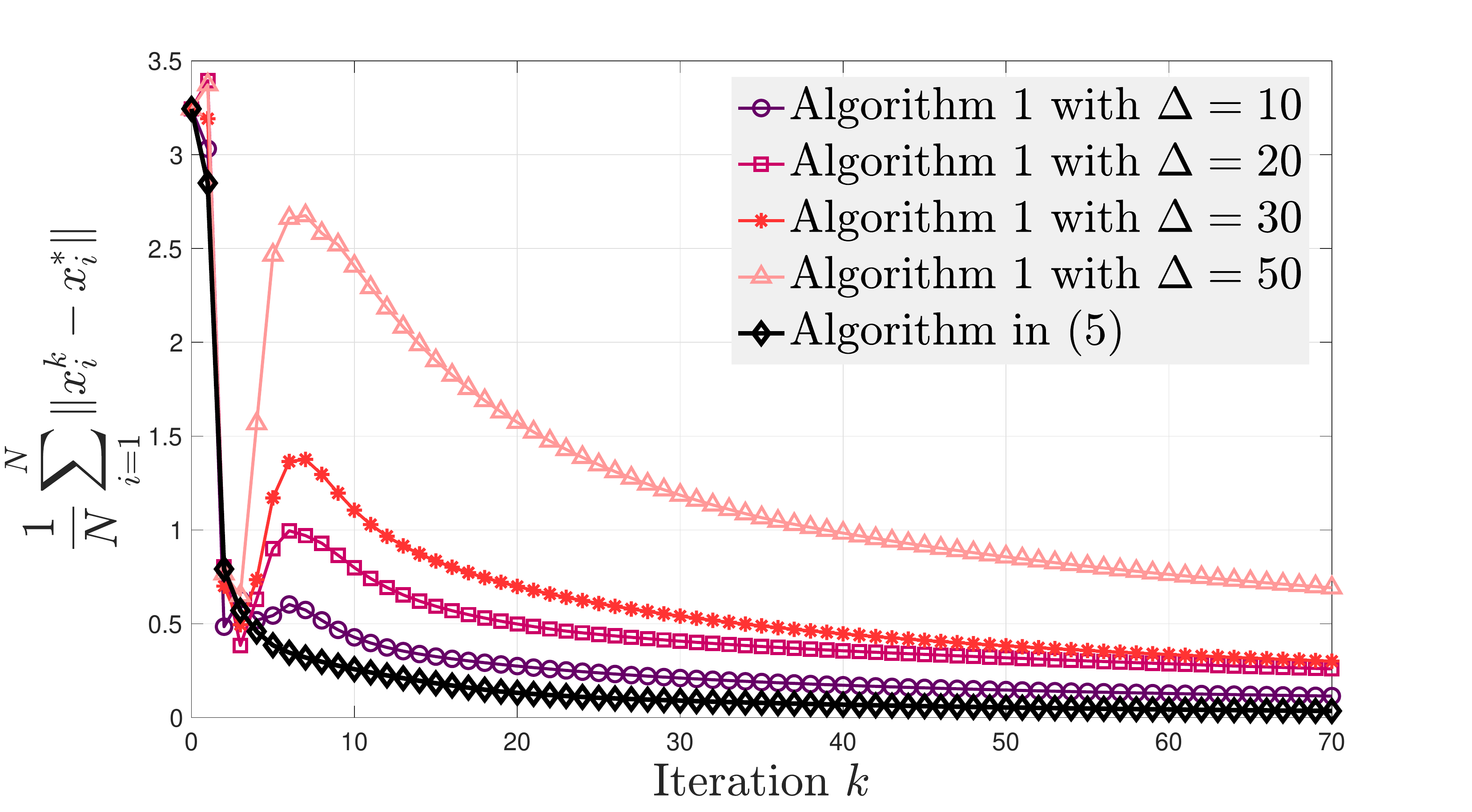}
    \caption{Iterates generated by Algorithm~\ref{Algo:PrivDNEComp} versus the Algorithm in \eqref{Eq:DAlgo1} for $\Delta = \{10, 20, 30, 50\}$.}
    \label{fig:xk_example}
\end{figure}

%\textit{Convergence Slowdown:} 
Consider a Cournot competition with $N=10$ players over $\Gfrak$ described in Figure \ref{fig:graph_example}. Player $i$'s cost is given by
$$c_i(x_i) = \zeta_{i,2} x_i^2 + \zeta_{i,1} x_i.$$ 
The cost coefficients are drawn randomly from
$$\zeta_{i,2} \sim {\sf unif}[0,1/2], \ \zeta_{i,1} \sim {\sf unif}[0, 1]$$ 
for each $i$. The strategy sets are identically $\Xcal_i = [0,5]$ for each $i$. Choose $\delta = \frac{1}{10}$ that parameterizes the matrix $W$. Let the price vary with demand $D$ as 
$$p(D) = 6 - \frac{1}{10}D.$$   

We initialize the algorithm with $x^0 = 1$ identically for all players. We use secure multi-party computing technique in \cite{gade2018acc} to design obfuscation sequence $\v{r}$ that satisfies \eqref{Eq:PrivateAlgo0} and 
 $$| r_{ij}^k | \leq \Delta.$$ 
The trajectory of the average distance of $x^k_i$'s from $x^*_i$ across players with $\alpha_k := (k+1)^{0.51}$ is shown in Figure \ref{fig:xk_example}.

Our algorithm converges to the equilibrium similar to the non-private algorithm in \eqref{Eq:DAlgo1}. However, its convergence is slower as seen in Figure~\ref{fig:xk_example}. The slowdown is especially pronounced for large $\Delta$'s and is an artifact of perturbations added by players to obfuscate information from the adversary. Thus, our algorithm design achieves privacy and asymptotic convergence to equilibrium, but sacrifices speed of convergence. An analytical characterization of the slowdown defines an interesting direction for future work.

\section{Proof of Theorem \ref{Th:Main}} \label{Sec:proofs}
%In Section \ref{Sec:PrivacyProof}, we show that our algorithm is private. In Section \ref{Sec:ConvergenceProof}, we sketch the proof that the iterates of the algorithm converge to the Nash equilibrium asymptotically. Details are omitted due to space limitations.

%----

\subsection{Proving Algorithm \ref{Algo:PrivDNEComp} is Private} \label{Sec:PrivacyProof}

Recall that $\Gfrak(\Acal^c)$ is the graph over non-adversarial nodes $\Acal^c$. Suppose $\Acal^c$ has $M$ nodes.
Let $I, J$ be two players in $\Acal^c$ and
$$F:={(f_i, \Xcal_i)}_{i \in \Vcal}, \quad \tilde{F} :={(\tilde{f}_i, \tilde{\Xcal}_i)}_{i \in \Vcal},$$
be two games in $\Fcal$ such that $\tilde{F}$ is identical to $F$, except that costs and strategy sets of players $I$ and $J$ are switched: 
$$ \tilde{f}_I = f_J, \ \tilde{f}_J = f_I, \ \tilde{\Xcal}_I = \Xcal_J, \ \tilde{\Xcal}_J = \Xcal_I.$$
For convenience, define $\pi: \Vcal \to \Vcal$ as the permutation that encodes the switch, i.e., 
$$\pi(I) = J, \ \pi(J)=I, \text{ and } \pi(i)=i \text{ for all } i \neq I,J.$$ 

Consider the execution of Algorithm \ref{Algo:PrivDNEComp} on $F$, given by
\begin{align*}
{\Esf}(F, \v{r}, \mathpzc{x}) := \{(x^k_i, v^k_i, \hat{v}^k_i) \text{ for } i \in \Vcal, k\geq 0 \},
\end{align*} 
with obfuscation sequence $\v{r}$ used in \eqref{Eq:PrivateAlgo1}, initialized with $\mathpzc{x} \in \cap_{i=1}^N \Xcal_i$. We prove that there exists an obfuscation sequence $\v{\tilde{r}}$ such that execution ${\Esf}(\tilde{F}, \v{\tilde{r}}, \mathpzc{x})$ of Algorithm \ref{Algo:PrivDNEComp} on $\tilde{F}$ with $\v{\tilde{r}}$ starting from $\mathpzc{x}$, is identical to ${\Esf}(F, \v{r}, \mathpzc{x})$, from ${\Asf}$'s perspective. 
An arbitrary permutation over $\Acal^c$ is equivalent to a composition of a sequence of switches among two players in $\Acal^c$. As a result, the algorithm execution on games in $\Fcal$ can be made to appear identical from ${\Asf}$'s standpoint, proving the privacy of Algorithm \ref{Algo:PrivDNEComp}.

In the rest of the proof, we show how to construct $\v{\tilde{r}}$ that ensures ${\Esf}(\tilde{F}, \v{\tilde{r}}, \mathpzc{x})$ and ${\Esf}(F, \v{r}, \mathpzc{x})$ appear identical to ${\Asf}$.

Adversary observes $\{x^k_j, v^k_j, \hat{v}^k_j\}$ for all $j\in\Acal$ at each $k\geq0$. Consequently, perturbations utilized by corrupted nodes $j\in\Acal$ are same in both executions,
\begin{align}
\tilde{r}^k_{ji} = r^k_{ji} \text{ for all } j\in\Acal.
\label{eq:q.r}
\end{align}

%\noindent 
Moreover, ${\Asf}$ observes $\hat{v}^k_j$ for all $j\in\Acal$ and hence, all messages received by $j\in\Acal$ from $i\in\Acal^c$, denoted by $v^k_{ij}$, are identical for both executions, i.e., 
\begin{align}
\begin{aligned}
 \tilde{v}^k_{ij} = {v}^k_{ij} 
& \iff \tilde{v}^k_i + \alpha^k \tilde{r}^k_{ij} = v^k_i + \alpha^k r^k_{ij}  \\
& \iff \alpha^k \tilde{r}^k_{ij} = v^k_i + \alpha^k r^k_{ij} - \tilde{v}^k_i. %= v^k_j + \alpha^k r^k_{ji} - v^k_{\pi(j)}
\end{aligned}
\label{eq:v.v}
\end{align}
Adversary observes $\bar{x}^k$ for each $k\geq 0$. Enforcing
$$\tilde{x}^k_i = x^k_{\pi(i)}, \ \tilde{v}^k_i = v^k_{\pi(i)}, \ \hat{\tilde{v}}^k_i = \hat{v}^k_{\pi(i)},$$
results in $\bar{\tilde{x}}^k=\bar{x}^k$. 
%and ${\Asf}$ is unable to differentiate between games $F$, $\tilde{F}$. 
We have 
\begin{align}
&\hat{\tilde{v}}^k_i = \hat{v}^k_{\pi(i)} \nonumber \\
& \iff \sum_{j\in\Ncal_i} W_{ij} (\tilde{v}^k_j +\alpha^k \tilde{r}^k_{ji}) = \sum_{j \in \Ncal_{\pi(i)}} W_{\pi(i)j} (v^k_j +\alpha^k r^k_{j \pi(i)}) \nonumber \\
& \iff \sum_{j\in\Ncal_i \cap \Acal^c} \tilde{r}^k_{ji} = \frac{1}{\alpha^k \delta} \sum_{j \in \Ncal_{\pi(i)}} W_{\pi(i)j} (v^k_j +\alpha^k r^k_{j\pi(i)})  \nonumber \\
&  \qquad \qquad  - \frac{1}{\alpha^k \delta}\sum_{j\in\Ncal_i} W_{ij} {v}^k_{\pi(j)}  - \sum_{j\in \Ncal_i \cap \Acal} \tilde{r}^k_{ji}. \label{Eq:Constraint3}%
\end{align}%

%\item 
The obfuscation used by each player $i\in\Acal^c$ is locally balanced, and hence, we have
\begin{align}
 \sum_{j \in \Ncal_i} \tilde{r}^k_{ij} = 0 
 \iff \sum_{j \in \Ncal_i \cap \Acal^c} \tilde{r}^k_{ij} = - \sum_{j \in \Ncal_i \cap \Acal} \tilde{r}^k_{ij}. \label{Eq:Constraint4}%
\end{align}%
Let $\gamma$ be a vector of $\tilde{r}^k_{ij}$'s for $i,j\in\Acal^c$. In the sequel, let $\bone$ denote a vector of ones of appropriate dimension. For graph $\Gfrak(\Acal^c)$, define its oriented incidence matrix $B$, adjacency matrix $A$, degree matrix $D$, and the normalized graph Laplacian matrix $L$ as 
\begin{align*}
&B_{ij} = \begin{cases}
1, & \text{if node }i \text{ is head of edge }j,    \\
-1, & \text{if node }i \text{ is tail of edge }j, \\
0, &\text{otherwise},
\end{cases} \\
&A_{ij} = 
\begin{cases} 
1, & \text{if } (i,j) \text{ is edge in } \Gfrak(\Acal^c), \\
0, & \text{otherwise},
\end{cases} \\
&D = \diag(A \bone), \text{ and } L = I - D^{-1/2} A D^{-1/2}.%
\end{align*}%
Using the notation $z_+ \coloneqq \max\{z, 0\}$ and $z_- \coloneqq z_+ - z$ for a scalar $z$, define $B_+$ and $B_-$ as the matrices obtained from $B$, applying the respective operator componentwise.
Then, \eqref{Eq:Constraint3} - \eqref{Eq:Constraint4} can be written as
\begin{align*}
\underbrace{\begin{pmatrix} B_- & B_+ \\ B_+ & B_-\end{pmatrix}}_{:=T} \gamma = \underbrace{\begin{pmatrix} \xi^1 \\ \xi^2\end{pmatrix}}_{:=\xi} ,
\end{align*}
%where, $T^1$ and $T^2$ are boolean matrices with $N-\tau$ rows each, defined as 
%\begin{align*}
%    T^1_{im} &= \begin{cases}
%    1,  & \text{if $m^{th}$ entry of }  \gamma \text{ is } \tilde{r}^k_{ti}, \text{ for any } t\in \Acal^c, \\
%    0,  & \text{otherwise, } 
%    \end{cases} \\
%    T^2_{im} &= \begin{cases}
%    1,  & \text{if $m^{th}$ entry of }  \gamma \text{ is } \tilde{r}^k_{it}, \text{ for any } t\in \Acal^c, \\
%    0,  & \text{otherwise. } 
%    \end{cases} 
%\end{align*}
We prove that, 
\begin{align}
\rank T = \rank (T \ | \ \xi) = 2M - 1,
\label{eq:rankCond}
\end{align}
to show that $T\gamma=\xi$ admits at least one solution.

Notice that 
 \begin{align}
 \begin{aligned}
( \bone^\T \ | \ -\bone^\T )
  \begin{pmatrix} B_- & B_+ \\ B_+ & B_-\end{pmatrix}
% & = 
%\left( \bone^\T [B_- - B_+] \ | \ \bone^\T [B_+ - B_-] \right) \\
  = 
 \left( - \bone^\T B \ | \ \bone^\T B \right) = 0,
 \end{aligned}
 \label{eq:T.not.indep}
 \end{align}
proving that rows of $T$ are not linearly independent. Next, we show that $\rank T \geq 2M-1$. To that end, we have
\begin{align}
\rank T
&= \rank (T T^\T) \nonumber \\
&= \rank \begin{pmatrix} 
B_- B_-^\T + B_+ B_+^\T &\; B_- B_+^\T + B_+ B_-^\T \\ 
B_+ B_-^\T + B_- B_+^\T &\; B_- B_-^\T + B_+ B_+^\T\end{pmatrix} \nonumber \\
&\stackrel{(a)}= \rank \begin{pmatrix} 
D & A \\ 
A & D \end{pmatrix} \nonumber \\
&= \rank D + \rank (D - A D^{-1}A) \nonumber \\
&= M +  \rank (I_M - D^{-1/2}A D^{-1}A D^{-1/2}) \nonumber \\
&= M +  \rank (I_M - (I_M-L)^2) \nonumber \\
&= M +  \rank (2L - L^2) \nonumber \\
&\stackrel{(b)}\geq M +  \underbrace{\rank L}_{=M-1} + \rank (2I_M - L) - M \nonumber \\
&\stackrel{(c)}= M - 1 + \rank(2 I_M - L) \nonumber \\
&= 2M-1, \label{eq:rankT}
\end{align}
where $I_M$ is the $M \times M$ identity matrix.
Here, (a) follows from the definition of $B_+, B_-, D, A$, (b) follows from Sylvester's nullity theorem (see \cite{horn2012matrix}) and the fact that the rank of graph Laplapcian for the connected graph $\Gfrak(\Acal^c)$ on $M$ nodes is $M-1$. Furthermore, since $\Gfrak(\Acal^c)$ is not bipartite, the eigenvalues of $L$ are strictly less than 2, according to Lemma 1.7 in \cite{chung1997spectral}. Therefore, we have $\rank(2 I - L) = M$ that implies (c). Thus, \eqref{eq:T.not.indep} and \eqref{eq:rankT} together yield $\rank  T = 2M-1$.

For the augmented matrix $(T \ | \ \xi)$, we have
\begin{align}
2 M - 1 \ \leq \ \rank( T \ | \ \xi ) \ \leq \ 2M. \label{Eq:Proof1} 
\end{align}
In the above relation, the inequality on the left follows from our earlier proof that $\rank T  = 2M-1$. The one on the right follows from the fact that the augmented matrix has $2M$ rows. We demonstrate that rows of $( T \ | \ \xi)$ are linearly dependent to conclude \eqref{eq:rankCond}. From \eqref{eq:T.not.indep}, we deduce
\begin{align*}
\left(\bone^\T \  | \ - \bone^\T \right) ( T \ | \ \xi) 
%&= \left( \bone^\T T^1 - \bone^\T T^2 \ | \ \bone^\T \xi^1 - \bone^\T \xi^2 \right)\\
&= \left( 0 \ | \ \bone^\T \xi^1 - \bone^\T \xi^2 \right).
\end{align*}
Now, we show $\bone^T \xi^1 - \bone^T \xi^2 = 0$ to conclude the proof. In the following, $| \Zcal |$ computes the cardinality of a set $\Zcal$.
\begin{small}
\begin{align*}
&\bone^T \xi^1 - \bone^T \xi^2   \\
&=  \frac{1}{\alpha^k \delta} \sum_{i\in\Acal^c} \sum_{j \in \Ncal_{\pi(i)}} W_{\pi(i)j} (v^k_j +\alpha^k r^k_{j\pi(i)})  \nonumber \\
&\quad  - \frac{1}{\alpha^k \delta}\sum_{i\in\Acal^c}\sum_{j\in\Ncal_i} W_{ij} {v}^k_{\pi(j)}  - \sum_{i\in\Acal^c} \sum_{j\in \Ncal_i \cap \Acal} \tilde{r}^k_{ji} + \sum_{i\in\Acal^c} \sum_{j \in \Ncal_i \cap \Acal} \tilde{r}^k_{ij}\\
%%%%%%%%%%%%%%%%%%%%%%%%%%%%%%%%%%%%%%%%%%%%%%%%%%%%%%%%
&{=}  \underbrace{\frac{1}{\alpha^k \delta} \sum_{i\in\Acal^c} \left[ \sum_{j \in \Ncal_{\pi(i)}} W_{\pi(i)j} v^k_j - \sum_{j\in\Ncal_i} W_{ij} {v}^k_{\pi(j)} \right]}_{:=Q^1}  \nonumber \\ %
&\quad + \underbrace{\sum_{i\in\Acal^c} \sum_{j \in \Ncal_{\pi(i)}}  r^k_{j\pi(i)} - \sum_{i\in\Acal^c} \sum_{j\in \Ncal_i \cap \Acal} \tilde{r}^k_{ji} + \sum_{i\in\Acal^c} \sum_{j \in \Ncal_i \cap \Acal} \tilde{r}^k_{ij}}_{:=Q^2} \\
&= Q^1 + Q^2,
%%%%%%%%%%%%%%%%%%%%%%%%%%%%%%%%%%%%%%%%%%%%%%%%%%%%%%%%
\end{align*}%
\end{small}%
where we have used $r^k_{ii}=0$ and $W_{ij} = \delta$ for $(i,j)\in\Ecal$. Utilizing $\pi(i)=i$,  for all $i \neq I,J$, simplify $Q^1$ as
\begin{small}
\begin{align*}
&\alpha^k Q^1 \\
&= \frac{1}{\delta} \sum_{i\in\Acal^c} \left[ \sum_{j \in \Ncal_{i}} W_{ij} v^k_j  - \sum_{j\in\Ncal_i} W_{ij} {v}^k_{\pi(j)} \right]\\
&= \frac{1}{\delta} \sum_{i\in\Acal^c} \left[ \sum_{j \in \Ncal_{i} \cap \Acal^c} W_{ij} v^k_j  - \sum_{j\in\Ncal_i \cap \Acal^c} W_{ij} {v}^k_{\pi(j)} \right]\\
&= \frac{1}{ \delta} \sum_{i\in\Acal^c} \left[\left(1-(|\Ncal_i|-1)\delta\right) v^k_i + \sum_{j\in\Ncal_i \cap \Acal^c \setminus \{i\}} \delta v^k_j\right]\\
&\quad - \frac{1}{ \delta} \sum_{i\in\Acal^c}\left[ \left(1-(|\Ncal_i|-1)\delta \right) v^k_{\pi(i)} + \sum_{j\in\Ncal_i \cap \Acal^c \setminus \{i\}} \delta v^k_{\pi(j)}\right] \\
&= \sum_{i\in\Acal^c} \left[ \left( \frac{1}{\delta} - |\Ncal_i| + 1 \right) \left(v^k_i - v^k_{\pi(i)} \right)  +  \sum_{j\in\Ncal_i \cap \Acal^c  \setminus \{i\}} \left( v^k_j- v^k_{\pi(j)} \right)\right] \\
&= \left( |\Ncal_I \cap \Acal^c \setminus \{I\} | + \frac{1}{\delta} - |\Ncal_I| + 1\right) \left(v^k_I - v^k_J \right) \\
&\quad + \left( |\Ncal_J \cap \Acal^c \setminus \{J\} | + \frac{1}{\delta} - |\Ncal_J| + 1\right) \left(v^k_J - v^k_I \right)\\
&= \left( |\Ncal_I \cap \Acal^c \setminus \{I\} | -|\Ncal_I | \right) \left(v^k_I - v^k_J \right) \\
&\quad + \left(|\Ncal_J \cap \Acal^c \setminus \{J\}|   -|\Ncal_J| \right) \left(v^k_J - v^k_I \right).
\end{align*}%
\end{small}%

Next, simplify $Q^2$ as
\begin{small}
\begin{align*}
Q^2 
%%%%%%%%%%%%%%%%%%%%%%%%%%%%%%%%%%%%%%%%%%%%%%%%%%%%%%%%
&{=}  \sum_{i\in\Acal^c} \sum_{j \in \Ncal_{i}}  r^k_{ji} - \sum_{i\in\Acal^c} \sum_{j\in \Ncal_i \cap \Acal} \tilde{r}^k_{ji} + \sum_{i\in\Acal^c} \sum_{j \in \Ncal_i \cap \Acal} \tilde{r}^k_{ij} \\
%%%%%%%%%%%%%%%%%%%%%%%%%%%%%%%%%%%%%%%%%%%%%%%%%%%%%%%%
&\stackrel{(a)}{=}  \sum_{i\in\Acal^c} \sum_{j \in \Ncal_i\cap\Acal^c} r^k_{ji} + \sum_{i\in\Acal^c} \sum_{j \in \Ncal_i \cap \Acal} \tilde{r}^k_{ij}\\ 
%%%%%%%%%%%%%%%%%%%%%%%%%%%%%%%%%%%%%%%%%%%%%%%%%%%%%%%
&\stackrel{(b)}{=} 
\sum_{i\in\Acal^c} \sum_{j \in \Ncal_{i}\cap\Acal^c} r^k_{ji} + \sum_{i\in\Acal^c} \sum_{j\in\Ncal_i\cap\Acal} \left[\frac{1}{\alpha^k}(v^k_i - v^k_{\pi(i)}) + r^k_{ij} \right] 
\end{align*}
\end{small}
\begin{small}
\begin{align*}
%%%%%%%%%%%%%%%%%%%%%%%%%%%%%%%%%%%%%%%%%%%%%%%%%%%%%%%%
&\stackrel{(c)}{=}  
\underbrace{\sum_{i\in\Acal^c} \left( \sum_{j\in\Ncal_i\cap\Acal^c} r^k_{ji} +\sum_{j\in\Ncal_i\cap\Acal} r^k_{ij}\right)}_{=0} \\ 
&\quad + {\sum_{i\in\Acal^c} \sum_{j\in\Ncal_i\cap\Acal} \left[\frac{1}{\alpha^k}(v^k_i - v^k_{\pi(i)}) \right]} \\
%%%%%%%%%%%%%%%%%%%%%%%%%%%%%%%%%%%%%%%%%%%%%%%%%%%%%%%
%%%%%%%%%%%%%%%%%%%%%%%%%%%%%%%%%%%%%%%%%%%%%%%%%%%%%%%
&\stackrel{(d)}{=} \frac{1}{\alpha^k}\sum_{i\in\Acal^c} |\Ncal_i\cap\Acal| (v^k_i - v^k_{\pi(i)}) \\
&= \frac{1}{\alpha^k}|\Ncal_I \cap \Acal| \left(v^k_I - v^k_J \right) + \frac{1}{\alpha^k}|\Ncal_J \cap \Acal| \left(v^k_J - v^k_I \right).
\end{align*}%
\end{small}%
Here, (a) follows from $r^k_{ji} = \tilde{r}^k_{ji}$ for all $j\in\Acal$ from \eqref{eq:q.r}. The equality in (b) follows from \eqref{eq:v.v}, (c) from \eqref{Eq:Constraint4}, and (d) from the properties of permutation $\pi$.
Combining the expressions for $Q^1$ and $Q^2$, we get
\begin{align*}
&Q^1 + Q^2 \\
&= \frac{1}{\alpha^k}\left[ |\Ncal_I \cap \Acal^c \setminus \{I\} | -|\Ncal_I \right] \left(v^k_I - v^k_J \right) \\
&\quad + \frac{1}{\alpha^k}\left[|\Ncal_J \cap \Acal^c \setminus \{J\}|   -|\Ncal_J| \right]. \left(v^k_J - v^k_I \right) \\
&\quad + \frac{1}{\alpha^k}|\Ncal_I \cap \Acal| \left(v^k_I - v^k_J\right) + \frac{1}{\alpha^k}|\Ncal_J \cap \Acal| \left(v^k_J - v^k_I \right) \\
&= 0,
\end{align*}
where the last line leverages the relation 
$$|\Ncal_i \cap \Acal| + |\Ncal_i \cap \Acal^c \setminus \{i\} | = |\Ncal_i|-1$$ for $i = I, J$. This completes the proof of privacy of our algorithm. $\hfill \qed$

\subsection{Proving Algorithm \ref{Algo:PrivDNEComp} Converges to Nash Equilibrium} \label{Sec:ConvergenceProof}
The non-expansiveness of the projection operator yields
\begin{small}
\begin{align}
    &\|x^{k+1}_i - x^*_i\|^2  \nonumber\\
    &= \|\proj_{\mathcal{X}_i}[x^{k}_i - \alpha^{k} \nabla_{x_i} f_i(x^{k}_i, N \hat{v}^{k}_i)] - x^{*}_i\|^2 \nonumber\\
    &= \|\proj_{\mathcal{X}_i}[x^{k}_i - \alpha^{k} \nabla_{x_i} f_i(x^{k}_i, N \hat{v}^{k}_i)]  \nonumber\\
    & \qquad \qquad - \proj_{\mathcal{X}_i}[x^{*}_i - \alpha^{k} \nabla_{x_i} f_i(x^{*}_i, \bar{x}^*)]\|^2 \nonumber\\
    &\leq \|x^k_i - x^*_i - \alpha^k (\nabla_{x_i} f_i(x^{k}_i, N \hat{v}^{k}_i) - \nabla_{x_i} f_i(x^{*}_i, \bar{x}^*))\|^2 \nonumber\\
 	&= \|x^k_i - x^*_i\|^2 + \underbrace{[\alpha^k]^2 \|\nabla_{x_i} f_i(x^{k}_i, N \hat{v}^{k}_i) - \nabla_{x_i} f_i(x^{*}_i, \bar{x}^*)\|^2}_{:=T^1_i}\nonumber\\
    & - \underbrace{2 \alpha^k (\nabla_{x_i} f_i(x^{k}_i, N \hat{v}^{k}_i) - \nabla_{x_i} f_i(x^{*}_i, \bar{x}^*))^\T(x^k_i - x^*_i).}_{:=T^2_i} \label{eq:x.ineq}%
\end{align}%
\end{small}%
Owing to the compactness of $\Xcal$'s, gradients $\nabla_{x_i} f_i$ are bounded. Such a bound, together with triangle inequality, yields an upper bound on $T^1_i$ as 
\begin{align}
T^1_i \leq [\alpha^k]^2 C^2.
\label{eq:T1.bound}
\end{align}
Define 
$$y^k := \frac{1}{N} \sum_{i=1}^N v^k_i, \quad C' := \max_{i} \max_{x_i \in \Xcal_i} \|x^k_i - x^*_i\| $$ 
and bound $T^2_i$ as 
\begin{small}
\begin{align}
T^2_i &= 2 \alpha^k \left[\nabla_{x_i} f_i(x^{k}_i, N \hat{v}^{k}_i) - \nabla_{x_i} f_i(x^{k}_i, N y^k)\right]^\T(x^k_i - x^*_i) \nonumber \\
   & \quad + 2 \alpha^k \left[\nabla_{x_i} f_i(x^{k}_i, N y^{k}) - \nabla_{x_i} f_i(x^{*}_i, \bar{x}^*)\right]^\T(x^k_i - x^*_i) \nonumber  \\
   &\geq -2\alpha^k N \bar{L} C' \|\hat{v}^k_i - y^k\|  \nonumber  \\
   & \quad + 2 \alpha^k \left[\nabla_{x_i} f_i(x^{k}_i, N y^{k}) - \nabla_{x_i} f_i(x^{*}_i, \bar{x}^*)\right]^\T(x^k_i - x^*_i),
\label{eq:T2.bound}%
\end{align}%
\end{small}%
where we use Cauchy-Schwarz inequality and Lipschitz continuity of $\nabla_{x_i} f_i$. 
To further simplify the bounds on $T^2_i$, we show that $Ny^k = \bar{x}^k$ using induction as follows. For $k = 0$, the relation follows from $v^0_i = x^0_i$. Assume that it holds for $k = 1, \ldots, K$, i.e., $Ny^K = \bar{x}^K$. Then, we have
\begin{small}
\begin{align}
N y^{K+1}
&= \sum_{i=1}^N v_i^{K+1} \nonumber\\
&\stackrel{(a)}= \sum_{i=1}^N  \left( \hat{v}_i^{K} + x_i^{K+1} - x_i^{K} \right) \nonumber\\
&\stackrel{(b)}= \sum_{i=1}^N  \left[ \sum_{j=1}^N W_{ij} \left(  v_j^K + \alpha^K r_{ji}^K \right) + x_i^{K+1} - x_i^{K} \right] \nonumber\\
&\stackrel{(c)}= \sum_{j=1}^N \underbrace{\sum_{i=1}^N W_{ij}}_{=1} v_j^K + \delta \alpha^K  \sum_{j=1}^N \underbrace{\sum_{i \in \Ncal_j}^N  r_{ji}^K}_{=0}  + \ \bar{x}^{K+1} - \bar{x}^{K} \nonumber\\
&\stackrel{(d)}= N y^K + \bar{x}^{K+1} - \bar{x}^{K} \nonumber\\
&= \bar{x}^{K+1},\label{Eq:ytracksxavg}%
\end{align}%
\end{small}%
where, (a) follows from \eqref{Eq:PrivateAlgo4}, (b) from  \eqref{Eq:PrivateAlgo1}, (c) from the doubly stochastic nature of $W$ and \eqref{Eq:PrivateAlgo0}. Finally, (d) follows from the induction hypothesis. 

Substitute $N y^{k} = \bar{x}^k$ in \eqref{eq:T1.bound} and combine that with \eqref{eq:T2.bound} in \eqref{eq:x.ineq}. The result, summed over $i \in \Vcal$ gives
\begin{small}
\begin{align}
    & \|x^{k+1} - x^*\|^2 \nonumber \\
    &\ \leq \|x^k - x^*\|^2 + [\alpha^k]^2 N C^2 + 2 \alpha^k N \bar{L}  C' \sum_{i=1}^N \|y^k - \hat{v}^k_i\| \nonumber \\
    & \quad \  - 2 {\sum_{i=1}^N \alpha^k \left[\nabla_{x_i} f_i(x^{k}_i, \bar{x}^{k}) - \nabla_{x_i} f_i(x^{*}_i, \bar{x}^*)\right]^\T(x^k_i - x^*_i)} \nonumber \\
    & \ = \|x^k - x^*\|^2 + [\alpha^k]^2 N C^2 + 2 \alpha^k N \bar{L} C' \sum_{i=1}^N \|y^k - \hat{v}^k_i\| \nonumber \\
    & \quad  \ - 2 \alpha^k \left[ \phi(x^k) - \phi(x^*) \right]^\T (x^k - x^*).
    \label{Eq:DissipationInequality}%
\end{align}%
\end{small}%	

We bound one of the terms on the right-hand side of the above relation in the next result.
\begin{lemma}
$\sum_{k=0}^\infty \alpha^k \|y^k - \hat{v}^k_i\| < \infty$, for all $i\in\Vcal$. \label{Lem:errors.converge}
\end{lemma}
The proof relies on the doubly stochastic nature of $W$ and two properties of obfuscation sequence -- boundedness of $\v{r}$ and balancedness property from \eqref{Eq:PrivateAlgo0}. We omit the proof due to space limitations. 
%We observe that \eqref{Eq:DissipationInequality} has the same structure as the dissipation inequality in Theorem~1, \cite{robbins1985convergence}. 
The square summability of $\alpha$'s, Lemma \ref{Lem:errors.converge} along with \eqref{Eq:DissipationInequality} allow us to infer that
 $\|x^k - x^*\|^2$ converges and 
$$\sum_{k=0}^\infty \alpha^k \Phi(x^k) \coloneqq \sum_{k=0}^\infty \alpha^k \left[ \phi(x^k) - \phi(x^*) \right]^\T (x^k - x^*) < \infty$$
using Theorem 1 in \cite{robbins1985convergence}.
The $\alpha$-sequence is nonsummable and $\phi$ is strictly monotone. Therefore, we have
$$ \liminf_{k\to \infty}\Phi(x^k) = 0.$$ 
The sequence of $x^k$'s remains bounded. Consider its bounded subsequence $x^{k^\ell}$ along which 
$$ \lim_{\ell \to \infty}\Phi(x^{k^\ell}) = \liminf_{k\to \infty}\Phi(x^k) = 0.$$
This subsequence admits a convergent subsequence, along which $\Phi$ goes to zero. Strict monotonicity of $\phi$ implies that this subsequence converges to $x^*$. Recall that $\|x^k - x^*\|^2$ converges, and this distance converges to zero over said subsequence, implying $\lim_{k\to \infty} x^k = x^*$. This completes the proof of Theorem \ref{Th:Main}. $\hfill \qed$

\section{Conclusions}
In this paper, we considered aggregate games played by agents that communicate over a network, each with private information. We showed that distributed algorithms for equilibrium computation in the literature are not designed with privacy requirements in mind, and consequently leak private information about players against honest-but-curious adversaries. Our proposed algorithm for NE computation exploits correlated perturbations to obfuscate aggregate estimates shared over the network. The algorithm asymptotically converges to the Nash Equilibrium. If the graph connecting non-adversarial players is connected and not bipartite, we show that our algorithm protects private information of non-adversarial players.

\bibliography{Central}

%\appendix
%\input{appendix_a.tex}
%

\end{document}